\documentclass{amsart}
\usepackage{amsmath,amssymb,amsfonts,enumerate,amsthm}
\usepackage[pdftitle={Title of Article}, verbose,
            colorlinks=true,
            naturalnames=true,
            linkcolor=blue,
            citecolor=red,
            pdfstartview=XYZ]{hyperref}

\newcommand{\Hom}{\mbox{Hom}\,}
\newcommand{\Ext}{\mbox{Ext}\,}
\newcommand{\Tor}{\mbox{Tor}\,}
\newcommand{\Spec}{\mbox{Spec}\,}

\newcommand{\depth}{\mbox{depth}\,}
\newcommand{\width}{\mbox{width}\,}
\renewcommand{\dim}{\mbox{dim}\,}

\newcommand{\pd}{\mbox{proj.dim}\,}
\newcommand{\id}{\mbox{inj.dim}\,}

\newcommand{\gid}{\mbox{Gid}\,}
\newcommand{\Gid}{\mbox{Gid}\,}

\newcommand{\lo}{\longrightarrow}

\newcommand{\fp}{\mathfrak{p}}

\newcommand{\fq}{\mathfrak{q}}

\newtheorem{thm}{Theorem}[section]
\newtheorem{cor}[thm]{Corollary}

\newtheorem{defn}[thm]{Definition}

\numberwithin{equation}{section}

\begin{document}
\bibliographystyle{amsplain}

\title{A Note On Gorenstein Injective Dimension}

\author{Leila Khatami}
\address{Leila Khatami\\ The Abdus Salam ICTP, Strada Costiera 11, 34100
Trieste, Italy}

\thanks{Email:  \href{mailto:lkhatami@ictp.it}{lkhatami@ictp.it}}



\begin{abstract}
Suppose that $M$ is a module over a commutative noetherian ring $R$.
It is proved that the Gorenstein injective dimension of $M$, if finite, 
equals $\sup 
\,  
\{\depth R_\fp - \width_{R_\fp}M_\fp \, | \, \fp \in \Spec(R)\}.$  

\end{abstract}

\maketitle

\section*{Introduction}

In 1976, Chouinard gave a general formula for injective dimension
of a module, when it is finite (cf. \cite{Chouinard}). \\
{\bf Theorem.} Let $M$ be an $R$-module of finite injective
dimension. Then $$\id_R M= \sup \{ \depth R_\fp -
\width_{R_\fp}M_{\fp}\, | \, \fp \in \Spec(R) \}.$$ Recall that
for a module $M$ over a local ring $R$, $\width_R M$ is defined as
the $\inf \{ \, i \, | \, \Tor_i^R(k,M) \neq 0\}$, where $k$ is
the residue field of $R$. This gives a general formula from which
the so called Bass formula can be concluded, namely, the injective
dimension of a finite module over a local ring is either infinite
or equals the depth of the base ring.

Our theorem \ref{GCH} extends the chouinard's formula for
Gorenstein injective dimension.
\\
{\bf Theorem.} Let $R$ be a noetherian ring and $M$ an $R$-module
of finite Gorenstein injective dimension. Then
$$\Gid_R M = \sup \{\depth R_{\fp}- \width_{R_\fp}M_\fp \, | \, \fp \in \Spec(R)\}.$$
\section{Main Theorem}
\begin{defn}
An $R$-module $G$ is said to be {\rm Gorenstein injective} if and
only if there exists an exact complex of injective $R$-modules,
$$I=\cdots \to I_2\lo I_1\lo I_0\lo I_{-1}\lo
I_{-2}\lo\cdots$$ such that the complex $\Hom_R(J,I)$ is exact for
every injective $R$-module $J$ and $G$ is the kernel in degree 0
of $I$. The {\rm Gorenstein injective dimension} of an $R$-module
$M$, $\Gid_R (M)$, is defined to be the infemum of integers $n$
such that there exists an exact sequence $$0 \to M \to G_0 \to
G_{-1} \to \cdots \to G_{-n} \to 0$$ with all $G_i$'s Gorenstein
injective.
\end{defn}




\begin{thm}\label{GCH}
Let $R$ be a commutative noetherian ring and $M$ an $R$-module of
finite Gorenstein injective dimension. Then $$\Gid_R(M)=\sup \{
\depth R_\fp-\width_{R_\fp}M_\fp \, | \, \fp \in \Spec(R) \}.$$
\end{thm}
\begin{proof}

First assume that $\Gid_R(M)=0$. By definition, there is an exact
sequence
$$E_\bullet: \, \, \cdots \to E_1 \to E_0 \to M \to 0$$ such that
every $E_i$ is injective. Set $K_i= \ker(E_{i-1} \to E_{i-2})$. \\
For any $\fp \in \Spec(R)$ and every $R_\fp$-module $T$, we have
$\Ext^i_{R_\fp}(T,M_\fp) \cong \Ext^{i+t}_{R_\fp}(T,(K_t)_\fp)$
for any two positive integers $i$ and $t$. Hence using
\cite[5.3(c)]{CFF} we get
$$\begin{array}{ll}
0 &= \sup \{ \, i \, |\, \Ext_{R_\fp}^i(T,M_\fp) \neq 0, \, \,
\mathrm{for} \, \, \mathrm{some} \, \, R_\fp \mathrm{-module} \,
\,
T \, \, \mathrm{with} \, \, \pd_{R_\fp} T < \infty \} \\
& \geq \depth R_\fp-\width_{R_\fp}M_\fp.
\end{array}$$
In addition, if $\fp$ is such that $\dim R/\fp=\dim_RM$ then,
using \cite[5.3(c)]{CFF} again, we get $\depth
R_\fp-\width_{R_\fp}M_\fp=0$.
\\
Therefore, $\sup \{\depth R_\fp-\width_{R_\fp}M_\fp \, | \, \fp
\in \Spec(R)\}=0$ for a Gorenstein injective module $M$.

Now assume that $n=\Gid_RM >0$. By \cite[2.16]{CFH}, there exists
a short exact sequence $$0 \to K \to L \to M \to 0,$$ where $K$ is
Gorenstein injective and $\id_RL=\Gid_RM=n$.
\\
Thus $$0=\sup\,  \{ \, \depth R_\fp-\width_{R_\fp}K_\fp \, | \,
\fp \in \Spec(R)\}.$$ On the other hand, by Chouinard's equality
\cite[3.1]{C}, we have
$$\id_RL =\sup \{\depth R_\fp-\width_{R_\fp}M_\fp \, | \, \fp \in \Spec(R)
\}.$$

For any $\fp \in \Spec(R)$, we denote $k(\fp)=R_\fp/{\fp R_\fp}$.
Then, for any $\fp \in \Spec(R)$, the exact sequence $0 \to K_\fp
\to L_\fp \to M_\fp \to 0$ induces the long exact sequence
$$\cdots \to \Tor_i^{R_\fp}(k(\fp),K_\fp) \to
\Tor_i^{R_\fp}(k(\fp),L_\fp) \to \Tor_i^{R_\fp}(k(\fp),M_\fp) \to
\Tor_{i-1}^{R_\fp}(k(\fp),K_\fp) \to \cdots$$ This sequence gives
rise to the following inequalities.
$$\begin{array}{c} \width_{R_\fp}L_\fp \geq
\min\{\width_{R_\fp}M_\fp,\width_{R_\fp}K_\fp\} \, \, \mathrm{and}\\
\width_{R_\fp}M_\fp \geq
\min\{\width_{R_\fp}L_\fp,\width_{R_\fp}K_\fp+1\}
\end{array}$$

If $\fp \in \Spec(R)$ is such that $\width_{R_\fp}K_\fp >
\width_{R_\fp}M_\fp$ then, by the mentioned inequalities,
$\width_{R_\fp}M_\fp = \width_{R_\fp}L_\fp.$
\\
But for any $\fp$ with $\width_{R_\fp}K_\fp \leq
\width_{R_\fp}M_\fp$, we also have $\width_{R_\fp}K_\fp \leq
\width_{R_\fp}L_\fp$. Thus
$$\begin{array}{c}
\depth R_\fp - \width_{R_\fp}L_\fp \leq \depth R_\fp -
\width_{R_\fp}K_\fp \leq 0 \, \, \mathrm{and}\\
\depth R_\fp - \width_{R_\fp}M_\fp \leq \depth R_\fp -
\width_{R_\fp}K_\fp \leq 0.
\end{array}$$

Therefore we have
$$\begin{array}{r}
\Gid_R(M)=\id_R(L)=\\
\sup\{\depth R_\fp-\width_{R_\fp}L_\fp \, | \, \fp \in \Spec(R)\}= \\
\sup\{\depth R_\fp-\width_{R_\fp}L_\fp \, | \, \fp \, \,
\mathrm{with}
 \, \, \width_{R_\fp}M_\fp < \width_{R_\fp}K_\fp\}= \\
\sup\{\depth R_\fp-\width_{R_\fp}M_\fp \, | \, \fp \, \,
\mathrm{with} \, \, \width_{R_\fp}M_\fp < \width_{R_\fp}K_\fp\}= \\
\sup\{\depth R_\fp-\width_{R_\fp}M_\fp \, | \, \fp \in \Spec(R)\}.
\end{array}$$
\end{proof}
\begin{cor}
Let $M$ be an $R$-module and $\fp \subseteq \fq$ prime ideals of
$R$. If $\gid_{R_\fp}M_\fp < \infty$ then $$\gid_{R_\fp}M_\fp \leq
\gid_{R_\fq}M_\fq. $$
\end{cor}

\providecommand{\bysame}{\leavevmode\hbox
to3em{\hrulefill}\thinspace}


\begin{thebibliography}{10}
\bibitem{Chouinard}
L.~G.~Chouinard, \emph{On finite weak and injective dimension},
Proc. Amer. Math. Soc. \textbf{60} (1976), 57--60.

\bibitem{CFF}
L.~W.~Christensen, H-B.~Foxby and A.~Frankild, \emph{Restricted
homological dimensions and Cohen-Macaulayness}, J. Algebra
\textbf{251} (2002) (1), 479--502.


\bibitem{CFH}
L.~W.~Christensen, A.~Frankild, and H.~Holm, \emph{On Gorenstein
projective, injective and flat dimensions-A functorial description
with applications}, To appear in J. Algebra.



\end{thebibliography}
\end{document}